\documentclass[11pt, a4paper]{article}

\usepackage{inputenc}
\usepackage{graphicx}
\usepackage{amssymb}
\usepackage{mathtools}
\usepackage{amsmath}
\usepackage{comment}
\usepackage{caption}
\usepackage{subcaption}
\usepackage{algorithm}
\usepackage{algorithmic}

\newcommand{\R}{\mathbb R}

\newtheorem{definition}{Definition}
\newtheorem{theorem}{Theorem}

\newtheorem{example}{Example}

\begin{document}
	\begin{center}
		{\large\bf{Multivariate approximation by polynomial and generalised rational functions.}}
	\end{center}
	
	\begin{center}{R. D\'iaz Mill\'an, V. Peiris, N. Sukhorukova and J. Ugon}
	\end{center}
	
	



\begin{abstract}
In this paper we develop an optimisation based approach to multivariate Chebyshev approximation on a finite grid. We consider two models: multivariate polynomial approximation and multivariate generalised rational approximation. In the second case the approximations are ratios of linear forms and  the basis functions are not limited to monomials. It is already known  that in the case of multivariate polynomial approximation on a finite grid the corresponding optimisation problems can be reduced to solving a linear programming problem, while the area of multivariate rational approximation is not so well understood.
In this paper we demonstrate that in the case of multivariate generalised rational approximation the corresponding optimisation problems are quasiconvex. This statement remains true even when the basis functions are not limited to monomials.  Then we apply a bisection method, which is a general method for quasiconvex optimisation. This method converges to an optimal solution with given precision. We demonstrate that the convex feasibility problems appearing in the bisection method can be solved using linear programming. Finally, we compare the deviation error and computational time for multivariate polynomial and generalised rational approximation with the same number of decision variables.
\end{abstract}

{\bf keywords}:
	 generalised rational approximation, Chebyshev approximation, quasiconvex functions, bisection method.

{\bf MSC2010}: 
 90C25, 
 90C26,
 90C90,              
 90C47,            
 65D15.              


\section{Introduction}\label{sec:introduction}

In this paper we study Chebyshev (uniform) approximation problem, where the approximations are multivariate polynomials or multivariate generalised rational functions. Before moving to multivariate approximation, we would like to highlight the most important results in univariate approximation, since this review also provides the motivation for the current study. 

Univariate polynomial approximation problems were studied for many decades~\cite{chebyshev,remez57}.  When, however,  nonsmooth or non-Lipschitz functions are subject to approximation, polynomial approximation is not very efficient. One of the ways to overcome this issue is to use piecewise polynomials or polynomial splines~\cite{NR,nurnberger,Schumaker1968,sukhorukovaoptimalityfixed}. In the case of polynomials and piecewise polynomials with fixed knots (points of switching from one poynomial to another), the optimisation problem is convex and therefore existing optimisation techniques can be used to tackle this problem. When the location of the knots is unknown (free knots), the problem becomes nonconvex and the corresponding optimisation problems are much harder~\cite{SukUgonCrouzeix2019,Meinardus1989,NR,nurnberger,Mulansky92,SukUgonTrans2017,su10}. In 1996, N\"urnberger~\cite{FreeKnotsOpenProblem96} identified this problem as a very hard and important open problem in approximation. 

Rational approximation (approximation by a ratio of polynomials) was considered as a promising alternative to the free knot spline approximation. In particular, there are theoretical studies~\cite{lorentz1996constructive,shevshevPopov1987}  that demonstrate that approximation by rational functions and free-knot spline approximations are ``closely related to each other''.   It should be noted, however, that~\cite{lorentz1996constructive,shevshevPopov1987} refer to non-uniform approximation.

Rational approximation was a very popular topic in the 1950s-70s~\cite{Achieser1965,Boehm1964,Meinardus1967rational,Ralston1965Reme,Rivlin1962} (just to name a few).  Most of the results were extended to the so called generalised  rational approximation in the sense of Cheney and Loeb~\cite{cheney1964generalized}, where the approximations are constructed as ratios  of linear forms and the corresponding basis functions are not limited to monomials. 

There are two main ways to approach rational approximation. One group of approaches~\cite{Trefethen2018} is dedicated to ``nearly optimal'' solutions, whose construction is based on Chebyshev polynomials. This approach is very efficient and very popular. The extension of this approach to non-monomial basis functions, however, remains open. The second approach is based on modern optimisation techniques, in particular, it uses the fact that the corresponding optimisation problems are quasiconvex and can be solved using a general quasiconvex optimisation methods. In particular, in~\cite{AMCPeirisSukhSharonUgon} the authors use a well-known bisection method (see \cite{SL}) to solve these problems. The advantage of this approach is that it can be extended to non-monomial basis functions. In \cite{diazsukhorugon}, the authors use a projection-type algorithm for solving these problems. In the current paper, we extend the first approach to multivariate functions.

The paper is organised as follows.
In section~\ref{sec:preliminaries} we provide the necessary background. 
In section~\ref{sec:optimisation} we formulate the optimisation problems appearing in multivariate polynomial and generalised rational approximation. In section~\ref{sec:bisection} we provide the details of bisection method and its application to multivariate generalised approximation. Section~\ref{sec:experiments} contains the results of numerical experiments. Section~\ref{sec:conc} we provide the conclusions and identify possible future research directions.

\section{Preliminaries and motivation}\label{sec:preliminaries}
\subsection{Polynomial approximation}

We start by introducing the following definitions and notation.

\begin{definition}
An exponent vector 
$${\bf e}=(e_1,\dots,e_l)\in \mathbb{R}^l,~e_i\in \mathbb{N},~i=1,\dots,l$$ for ${\bf x}\in\mathbb{R}^l$
defines a {\em monomial}
$${\bf x^e} = x_1^{e_1} x_2^{e_2} \dots x_l^{e_l}.$$ 
\end{definition}

\begin{definition}
The \emph{degree} of a monomial $x^e$ is the sum of the components of $e$:
$$deg({\bf x}^e)=\sum_{i=1}^{l}e_i.$$
\end{definition}

\begin{definition}
A product $c {\bf x^e},$ where $c\ne 0$ is called a \emph{term}. A \emph{multivariate polynomial} is a sum of a finite number of terms. The \emph{degree} of a polynomial is the largest degree of its composing monomials.
\end{definition}

Assume that the dimension of the parameter space of a polynomial of degree $m$ is $d$. Note that in the case of linear functions and univariate function (that is, $l=1$) $d=l+1$. If $l\geq 2$ and $m\geq 2$ then $d$ (the total number of possible monomials of degree at least $m$) is increasing very fast. 
\begin{example}
Let $l=2$ (variables $x$ and $y$) and $m=2$. Then the total number of all possible monomials of degree zero is one, of degree one is~$l=2$ ($x$ and $y$) and of degree two is three ($x^2$, $y^2$ and $xy$). 
\end{example}
\subsection{Generalised rational approximation}
In this paper, we consider multivariate generalised rational approximation in terms of Cheney and Loeb~\cite{cheney1964generalized}. Namely, the approximations are the ratios of linear forms
\begin{equation}\label{eq:genrational}
R_{n,m}(x) =  P(x)/Q(x) = \frac{\sum_{i=0}^{n} p_i g_i(x)}{\sum_{j=0}^{m} q_jh_j(x)},
\end{equation}
where $g_i(x)$, $i=1,\dots,n$ and $h_j(x)$, $j=1,\dots,m$ are basis functions and $x\in X\in\R^l$, $l>1$.

In the next section, we formulate the optimisation problems, appearing in multivariate polynomial and generalised rational approximation.

\section{Optimisation problems}\label{sec:optimisation}
\subsection{Multivariate polynomial approximation}\label{ssec:polynomial}

Suppose that a continuous function $f({\bf x}):\mathbb{R}^l\rightarrow\mathbb{R}$ is to be approximated on a compact set $X\in \R^l$ by a function
 \begin{equation}\label{eq:model_function}
 L({\bf A,x})=a_0+\sum_{i=1}^{n}a_ig_i({\bf x}),
 \end{equation}
 where $g_i({\bf x})$ are the basis functions and the polynomial coefficients ${\bf A} = (a_0,a_1,\dots,a_n)$ are the decision variables. At a point \({\bf x}\) the deviation between the function \(f\) and the approximation is:
 \begin{equation}
   \delta({\bf A,x}) = |f({\bf x}) - L({\bf A,x})|.
\end{equation}
   \label{eq:deviation}
 Then the uniform approximation error over the set \(X\) is
 \begin{equation}
   \label{eq:uniformdeviation}
\Psi({\bf A})=\sup_{{\bf x}\in X}\|f({\bf x})-a_0-\sum_{i=1}^{n}a_ig_i({\bf x})\|.
\end{equation}
and the approximation problem is
 \begin{equation}\label{eq:obj_fun_con}
   \mathrm{minimise~}\Psi({\bf A}) \mathrm{~subject~to~} {\bf A}\in \R^{n+1}.
 \end{equation}
Note that $\Psi({\bf A})$ is convex (supremum of affine functions).

It is easy to see that this problem is convex and can be formulated as the following (possibly semi-infinite) linear programming problem:
\begin{equation}\label{eq:LPpol}
{\rm minimise}~ z
\end{equation}
subject to
\begin{equation}\label{eq:LPpolcon}
f({\bf x})-a_0-\sum_{i=1}^{n}a_ig_i({\bf x})\leq z,~
a_0-\sum_{i=1}^{n}a_ig_i({\bf x})-f({\bf x})\leq z,~{\bf x}\in X
\end{equation}
Necessary and sufficient optimality conditions for multivariate polynomial approximation can be found in~\cite{rice1963tchebycheff} and a more geometrical formulation can be found in~\cite{Sukhorukova2018}.

\subsection{Multivariate generalised rational approximation}\label{ssec:rational}

In the case of multivariate generalised rational approximation, the optimisation problem is
\begin{equation}\label{eq:problem}
\min_{\bf A,B}\sup_{{\bf x}\in X}\left|f({\bf x})-\frac{{\bf A}^T{\bf G}({\bf x})}{{\bf B}^T{\bf H}({\bf x})}\right|,
\end{equation}
subject to
\begin{equation}\label{eq:positivity}
{\bf B}^T{\bf H}({\bf x})>0,~{\bf x}\in X,
\end{equation}
where 
${{\bf A}}=(a_0,a_1,\dots,a_n)^T\in\R^{n+1}, {{\bf B}}=(b_0,b_1,\dots,b_m)^T\in\R^{m+1}$ are our decision variables, and ${\bf G}({\bf x})=(g_0({\bf x}),\dots,g_n({\bf x}))^T$, ${\bf H}({\bf x})=(h_0({\bf x}),\dots,h_m({\bf x}))^T$, where $g_j({\bf x})$, $j=1,\dots,n$ and $h_i({\bf x})$, $i=1,\dots, m$ are known basis functions defined on $X$.  Therefore,  the approximations are ratios of linear combinations of basis functions. 

The following theorem is the extension of the results from~\cite{SL,AMCPeirisSukhSharonUgon} for the case of multivariate approximation.
\begin{theorem}\label{thm:quasiconvexity}
Function
\begin{equation}\label{eq:phirational}
\Phi({\bf A,B})=\sup_{{\bf x}\in X}\left|f({\bf x})-\frac{{\bf A}^T{\bf G}({\bf x})}{{\bf B}^T{\bf H}({\bf x})}\right|
\end{equation}
is quasiconvex on $X$.
\end{theorem}
{\bf Proof:}
It is shown in~\cite{SL} that ratios of linear forms are quasilinear. The supremum of quasilinear functions is quasiconvex. 

\hskip 300pt
$\square$

Problem~(\ref{eq:problem})-(\ref{eq:positivity}) can be also formulated as follows:
\begin{equation}\label{eq:problemLP_obj}
\min~z
\end{equation}
subject to
\begin{equation}\label{eq:problemLP_con}
f({\bf x})-\frac{{\bf A}^T{\bf G}({\bf x})}{{\bf B}^T{\bf H}({\bf x})}\leq z,~\frac{{\bf A}^T{\bf G}({\bf x})}{{\bf B}^T{\bf H}({\bf x})}-f({\bf x})\leq z,~{\bf x}\in X,
\end{equation}
\begin{equation}\label{eq:positivityLP}
{\bf B}^T{\bf H}({\bf x})>0,~{\bf x}\in X.
\end{equation}

Without loss of generality, the constraint~(\ref{eq:positivityLP}) can replaced by
\begin{equation}\label{eq:delta}
{\bf B}^T{\bf H}({\bf x})\geq\delta,~{\bf x}\in X,
\end{equation}  
where $\delta$ is an arbitrary positive number.

If $z$ is fixed, $X$ is a finite grid then the constraint set~(\ref{eq:problemLP_con})-(\ref{eq:delta}) becomes a set of solutions to a finite number of nonstrict linear inequalities.

Note that if $({\bf A}$, ${\bf B})$ is an optimal solution, then  $(\lambda{\bf A}$, $\lambda{\bf B})$ is also an optimal solution. Therefore, it is reasonable to apply a normalisation method to ensure uniqueness  of the optimal solution. In our numerical experiments we fix one of the coefficients of the denominator (namely, the coefficient with the monomial of degree zero). In general, this approach should be applied with care, since in some cases this may lead to a zero denominator.  

There are a number of methods for solving quasiconvex minimisation problems. One of them, the bisection method~\cite{SL}, is described in the next section. This method relies on the fact that the sublevel sets of quasiconvex functions are convex. This property is also considered as one of the equivalent definitions of quasiconvex functions.

\section{Bisection method for multivariate generalised rational approximation}\label{sec:bisection}
The bisection method is a simple, but efficient approach for minimising quasiconvex functions~\cite{SL}.  To initiate the bisection method, one needs to define the following parameters.
\begin{itemize}
\item The absolute precision for maximal deviation~$\varepsilon$. 
\item The upper bound $u=\frac{\max_{{\bf x}\in X}f({\bf x})-\min_{{\bf x}\in X}f({\bf x})}{2}$.
\item The lower bound $l=0$.
\end{itemize}
Note that in some cases, more accurate $u$ and $l$ are available. For example one can use the best polynomial approximation to define $u$.

Set $z=\frac{1}{2}(u+l)$ and check  if the set of constraints~(\ref{eq:problemLP_con})-(\ref{eq:positivityLP}) has a feasible point. If this set is empty, update the lower bound~$l=z$, otherwise update the upper bound~$u=z$.  Repeat this procedure while $u-l\geq \varepsilon$.

In general, checking the existence of feasible points may be a difficult task (convex feasibility problems). There are a number of efficient methods (\cite{BauschkeLewis, Shi-yaXu, Zaslavski,YangYang, Zhao} just to  name a few), but there are still open problems and potential research direction to improve the efficiency of these methods. The discussion of the details of these method in application to generalised rational multivariate approximation is out of scope of the current paper, but this is one of our future research directions. Below we review two possible approaches.

\subsection{Solving the convex feasibility subproblem}

\subsubsection{Linear Programming}

In the case of multivariate generalised rational approximation, this problem  can be reduced to solving a linear programming problem. Indeed, note that the denominator of the approximation is positive,  fix~$z$ and solve the following problem:
\begin{equation}\label{eq:problemLP_obj_ax}
\min~\tilde{u}
\end{equation}
subject to
\begin{equation}\label{eq:problemLP_con_ax1}
f({\bf x}){{\bf B}^T{\bf H}({\bf x})}-{{\bf A}^T{\bf G}({\bf x})}\leq z{{\bf B}^T{\bf H}({\bf x})}+\tilde{u},
\end{equation}
\begin{equation}\label{eq:problemLP_con_ax2}
{{\bf A}^T{\bf G}({\bf x})}-f({\bf x}){{\bf B}^T{\bf H}({\bf x})}\leq z{{\bf B}^T{\bf H}({\bf x})}+\tilde{u},~{\bf x}\in X,
\end{equation}
\begin{equation}\label{eq:positivityLPax}
{\bf B}^T{\bf H}({\bf x})\geq\delta,~{\bf x}\in X.
\end{equation}
If an optimal solution $\tilde{u}\leq 0$, the set~(\ref{eq:problemLP_con})-(\ref{eq:positivityLP}) has a feasible point, otherwise the set is empty.

If $X$ is a finite grid, then (\ref{eq:problemLP_obj_ax})-(\ref{eq:positivityLPax}) is a linear programming problem and can be solved efficiently at each step of the bisection method.

\subsubsection{Splitting Method}

Alternatively, one can use a splitting algorithm, such as the alternating projection, or Douglas-Rachford methods.

Consider the following sets:
\[
\begin{aligned}
C^+(z) &= \Big\{({\bf A,B}):f({\bf x})-\frac{{\bf A}^T{\bf G}({\bf x})}{{\bf B}^T{\bf H}({\bf x})}\leq z,{\bf B}^T{\bf H}({\bf x})\geq\delta, : {\bf x}\in X\Big\},\\
C^-(z) &= \Big\{({\bf A,B}):\frac{{\bf A}^T{\bf G}{({\bf x})}}{{\bf B}^T{\bf H}{({\bf x})}}-f({{\bf x})\leq z,{\bf B}^T{\bf H}({\bf x})}\geq\delta,~{\bf x}\in X\Big\},
\end{aligned}
\]

These sets are both convex and nonempty. Convexity is due to the fact that they are the intersection of half-planes (as illustrated by equations~\eqref{eq:problemLP_con_ax1} and~\eqref{eq:problemLP_con_ax2} with $u=0$). Consider ${\bf A}^+ = (\max_{{\bf x}\in X} f({\bf x}),0,\ldots,0)$ and ${\bf B}^{+} = (1,0,\ldots,0)$. Clearly $({\bf A}^+,{\bf B}^+)\in C^+(z)$ for any $z\geq 0$.  Similarly if ${\bf A}^- = (\min_{{\bf x}\in X} f({\bf x}),0,\ldots,0)$ and ${\bf B}^{-} = (1,0,\ldots,0)$. Clearly $({\bf A}^-,{\bf B}^-)\in C^-(z)$ for any $z\geq 0$.

Therefore the subproblem consists in finding the (possibly empty) intersection of two nonempty convex sets, written as: 

\begin{equation}\label{eq:feasibi}
  ({\bf A^*,B^*})\in  C^+(z)\cap C^-(z).
 \end{equation}
 
When $X$ is a finite grid, the problem becomes in a finite intersection of hyperplanes. Since the projection onto a hyperplane has a closed formula,  the problem is easily solvable.
   
   A difficult problem is when $X$ is an interval or an infinite grid. For this case, the projection onto $C^+(z)$ and $C^-(z)$ is not a trivial problem and motive for our future research.

%
%

\section{Numerical experiments}\label{sec:experiments}
In this section we demonstrate a few numerical aspects of the optimisation problems appearing in section~\ref{sec:optimisation}. To illustrate the relative merits of multivaiate generalised rational approximation over multivariate polynomial approximation, we compare these approximations when the number of decision variables in both cases are the same. The code used in this section to develop the numerical results was implemented in MATLAB~R2020a.

\begin{figure}
	\centering
	\includegraphics[width=0.9\linewidth]{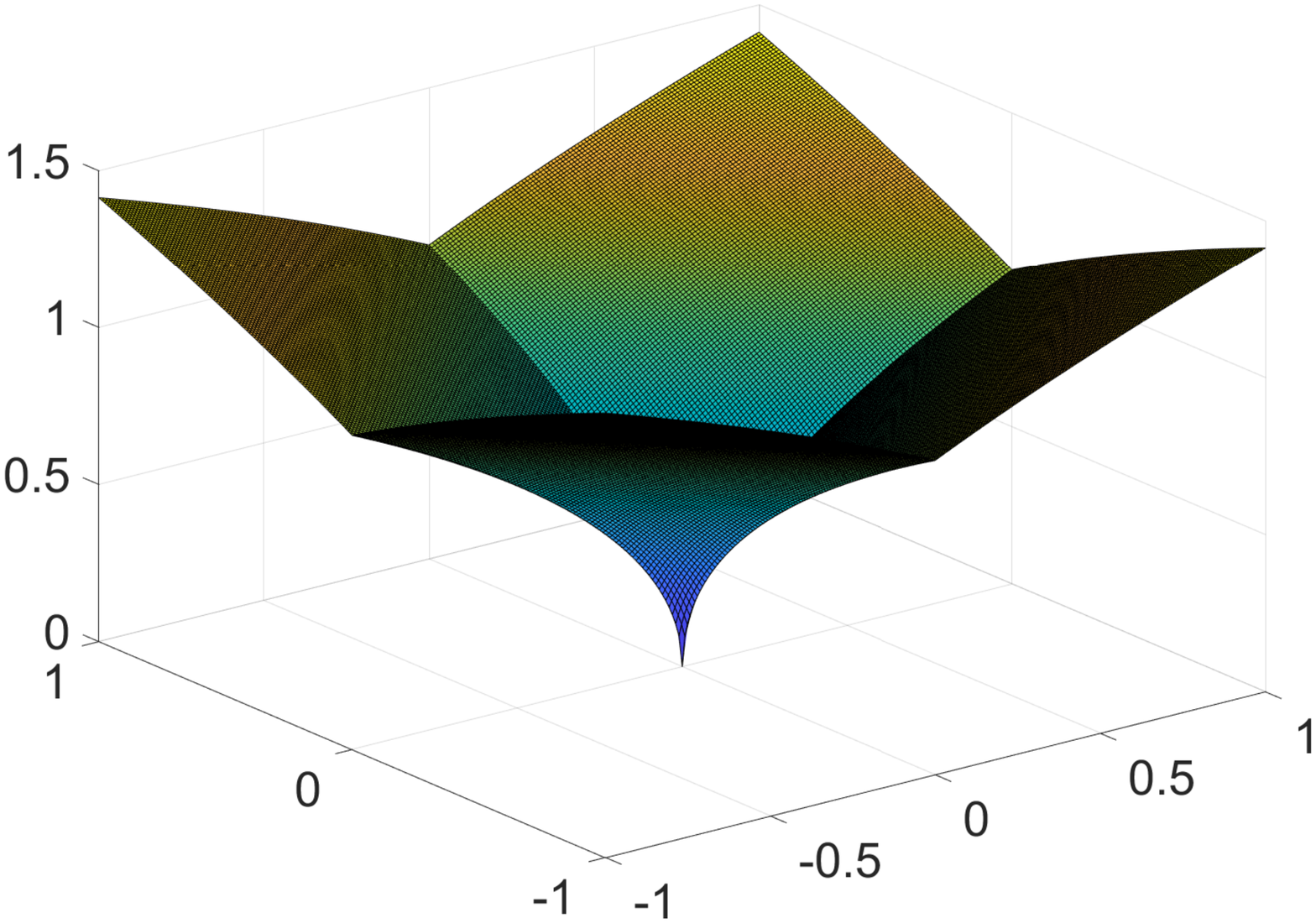}
	\caption{The original function $f(x,y)$.}
	\label{fig:original_function}
\end{figure}

We run our numerical experiments on a finite grid, the step size is~0.01, approximating 
\begin{equation*}
	f(x,y) = \sqrt{|x|+|y|}, \quad x,y \in [-1,1]
\end{equation*}
This is a nonsmooth and non-Lipschitz function with an abrupt change at the point $(0,0)$ (see figure~\ref{fig:original_function}).

\begin{figure}
	\begin{center}
		\begin{subfigure}{.49\textwidth}
			\includegraphics[width=\textwidth]{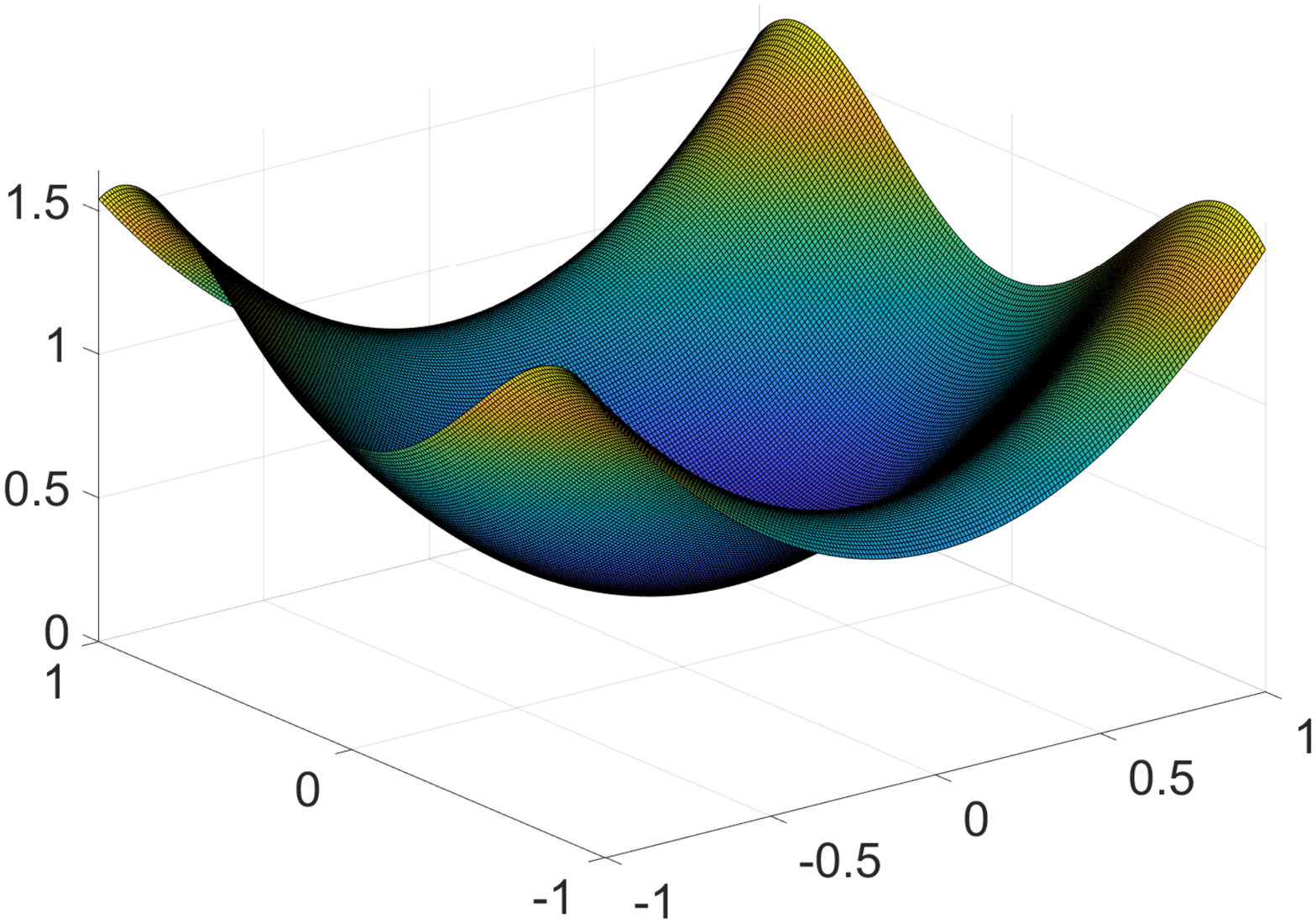}  
			\caption{Multivariate polynomial approximation.}
			\label{fig:polyapp}
		\end{subfigure}
		\begin{subfigure}{.49\textwidth}
			\includegraphics[width=\textwidth]{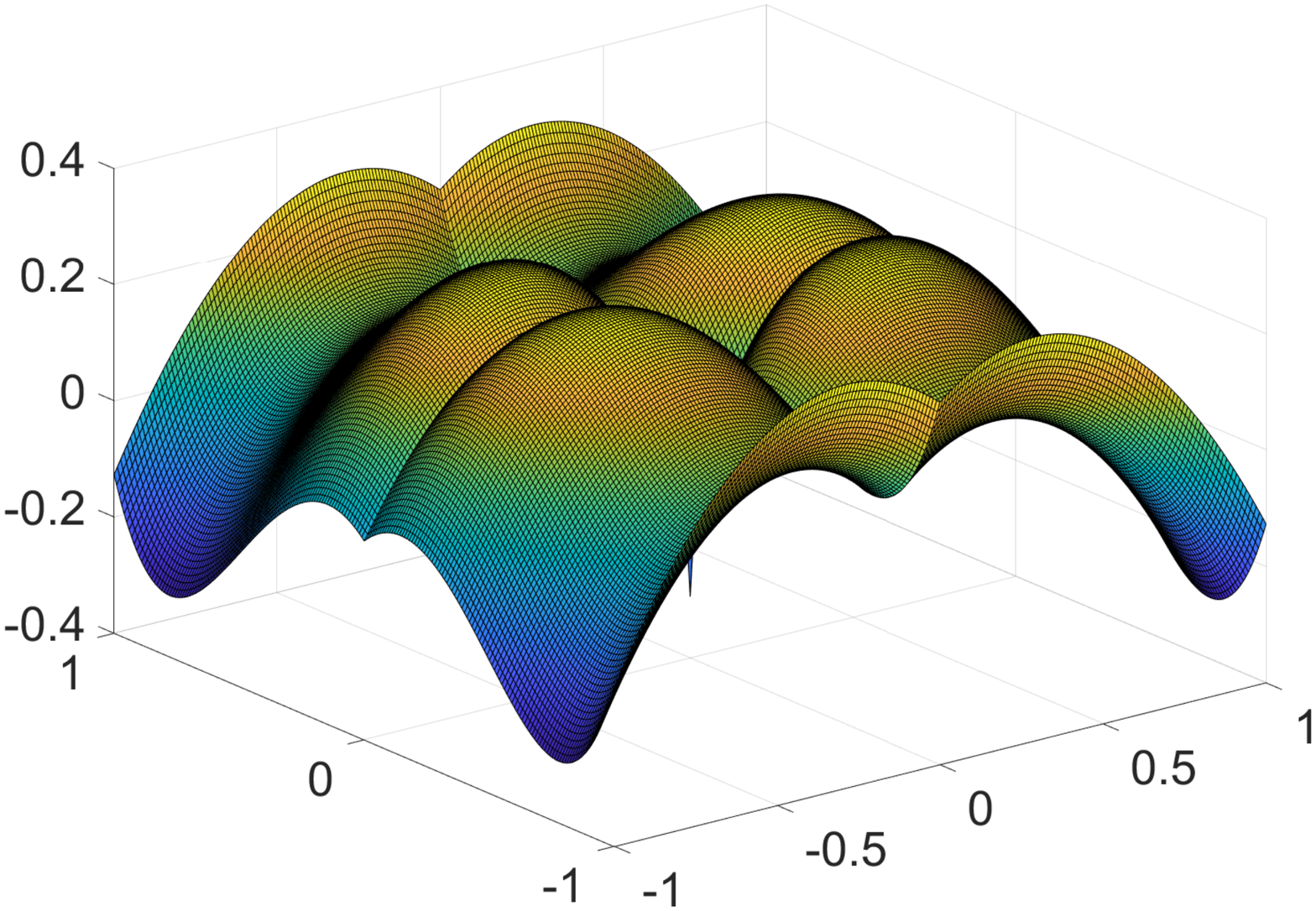}  
			\caption{Error curve.}
			\label{fig:polyerr}
		\end{subfigure}
		\caption{Multivariate polynomial approximation $P({\bf A},x,y)$ of degree 4 and its error curve.}
		\label{fig:multipolynomial}
	\end{center}
\end{figure}

Suppose that the function $f(x,y)$ is approximated by a multivariate polynomial of degree 4. The corresponding optimisation problem can be stated as follows:
\begin{equation*}
	\min_{\bf A} \max_{x,y \in [-1,1]} \left| f(x,y) - P({\bf A},x,y) \right|,
\end{equation*}
where
\begin{equation*}
	P({\bf A},x,y) = a_0 + a_1 x + a_2 y + a_3 x^2 + a_4 y^2 + a_5 xy + a_6 x^3 + a_7 y^3 + a_8 x^2y + a_9 xy^2 + a_{10} x^4
\end{equation*}
Vector ${\bf A}$ contains 11~components, that are also the decision variables. The result is presented in figure~\ref{fig:multipolynomial}, where figure~\ref{fig:polyapp} shows the multivariate polynomial approximation for the function $f(x,y)$ side by side with its maximum alternating error curve depicted in figure~\ref{fig:polyerr}.


The obtained coefficients for the multivariate polynomial $P(A,x,y)$ is presented in table~\ref{table:polyapp}. Only the basis functions of even degree correspond to nonzero coefficient. This is not surprising, since $f(x,y)$ is symmetric. Therefore, the resulting approximation can be written as,
$$P({\bf A},x,y) = 0.2944 + 1.6439 x^2 + 0.7638 y^2 - 1.1611 x^4.$$
The maximum deviation for the corresponding problem is 0.2944 and the computational time is 53.2848 seconds.

\begin{table}
	\centering
	\caption{Obtained values for the coefficients for $P({\bf A},x,y)$.}
	{\begin{tabular}{ |c|c|c| } 
		\hline
		Basis function & Coefficient & Value  \\ \hline 
		$1$ & $a_0$ & 0.2944 \\ \hline
		$x$ & $a_1$ & 0 \\ \hline
		$y$ & $a_2$ & 0  \\ \hline
		$x^2$ & $a_3$ & 1.6439  \\ \hline
		$y^2$ & $a_4$ & 0.7638  \\  \hline
		$xy$ & $a_5$ & 0 \\  \hline
		$x^3$ & $a_6$ & 0  \\  \hline
		$y^3$ & $a_7$ & 0 \\  \hline
		$x^2y$ & $a_8$ & 0 \\  \hline
		$xy^2$ & $a_9$ & 0 \\  \hline
		$x^4$ & $a_{10}$ &  -1.1611 \\  
		\hline
	\end{tabular}}
\label{table:polyapp}
\end{table}

\begin{figure}
	\begin{center}
		\begin{subfigure}{.49\textwidth}
			\includegraphics[width=\textwidth]{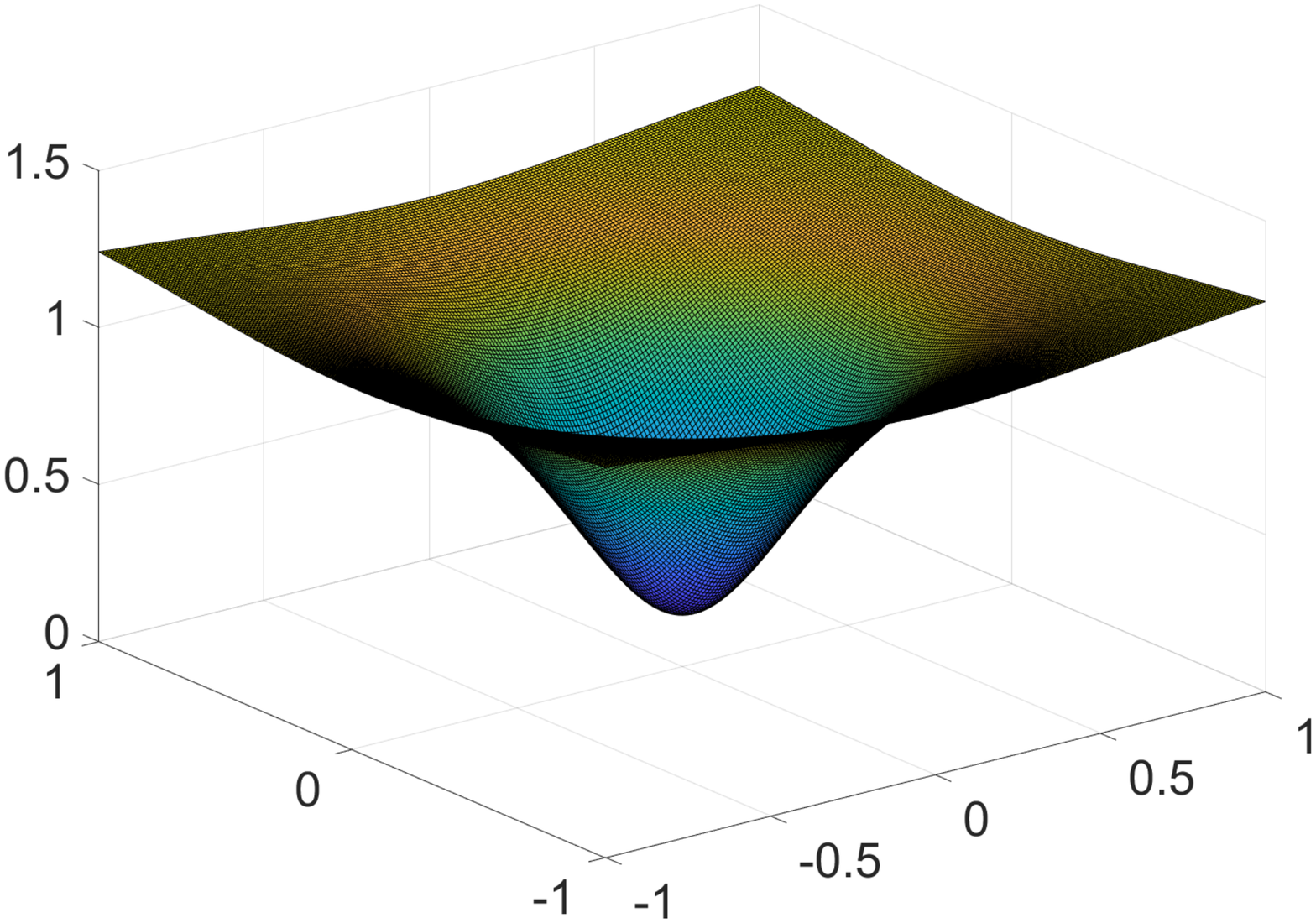}  
			\caption{Multivariate generalised rational approximation.}
			\label{fig:ratapp}
		\end{subfigure}
		\begin{subfigure}{.49\textwidth}
			\includegraphics[width=\textwidth]{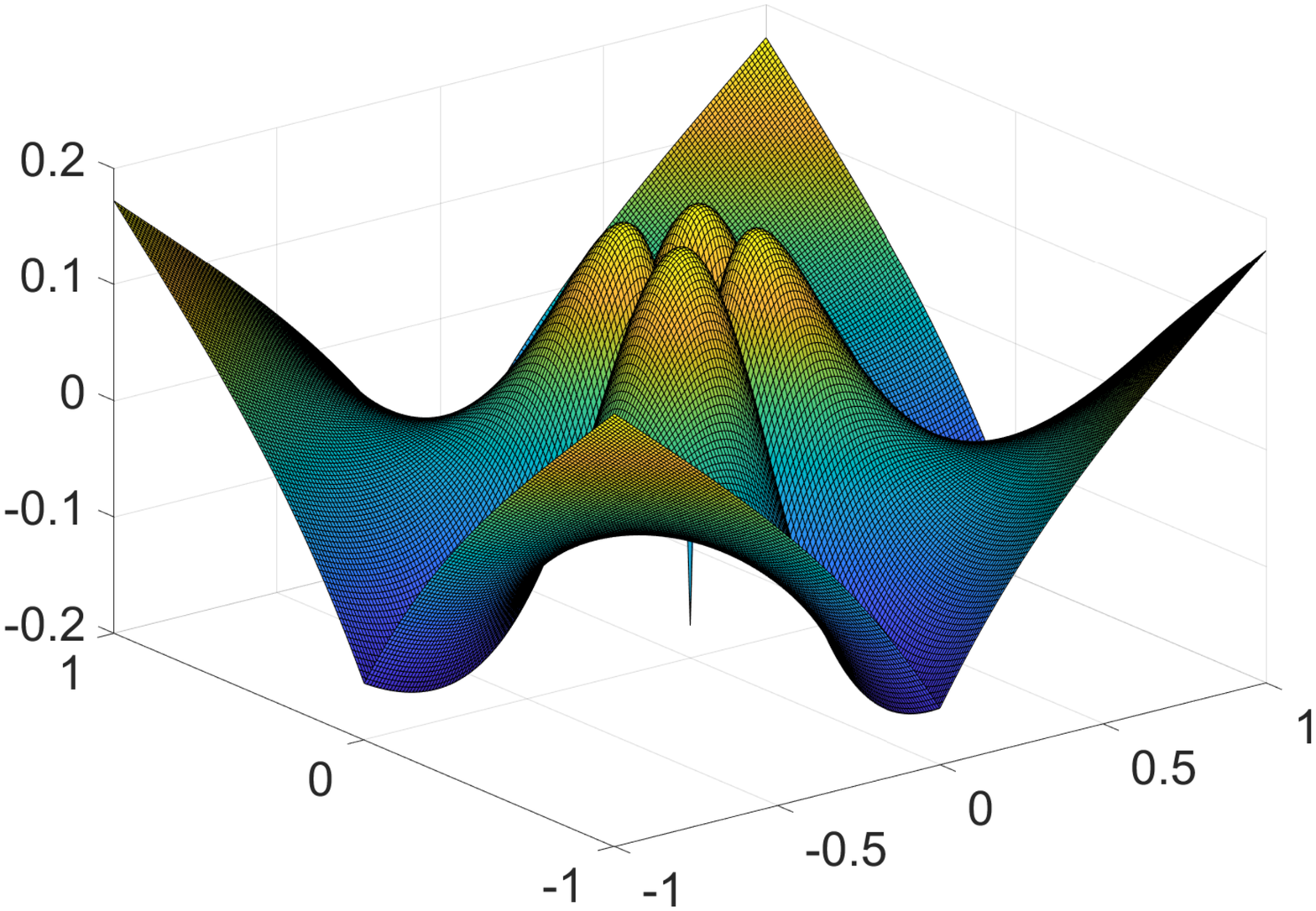}  
			\caption{Error curve.}
			\label{fig:raterr}
		\end{subfigure}
		\caption{Multivariate generalised rational approximation $P({\bf A},x,y)$ of degree 4 and its error curve.}
		\label{fig:multirational}
	\end{center}
\end{figure}

Now let us approximate the same function $f(x,y)$ by a multivariate generalised rational function of degree~2 in the numerator and degree~2 in the denominator. The corresponding optimisation problem is as follows:
\begin{equation*}
	\min_{{\bf A},{\bf B}} \max_{x,y \in [-1,1]} \left| f(x,y) - R({\bf A},{\bf B},x,y) \right|
\end{equation*}

\begin{equation*}
	R({\bf A},{\bf B},x,y) = \frac{a_0 + a_1 x + a_2 y + a_3 x^2 + a_4 y^2 + a_5 xy}{1 + b_1 x + b_2 y + b_3 x^2 + b_4 y^2 + b_5 xy}
\end{equation*}
This function $R({\bf A},{\bf B},x,y)$ has exactly 11~decision variables as in the multivariate polynomial case. We use the bisection method introduced in section~\ref{sec:bisection} to find the approximation. Figure~\ref{fig:ratapp} shows the multivariate generalised rational approximation for the function $f(x,y)$ whereas figure~\ref{fig:raterr} illustrates its maximum alternating error curve. By comparing figure~\ref{fig:multipolynomial} and figure~\ref{fig:multirational}, one can notice that the multivariate generalised rational function introduces a superior result over the multivariate polynomial one even though the dimension for both cases are the same.


The obtained coefficients for $R({\bf A},{\bf A},x,y)$ is presented in table~\ref{table:ratapp}. Hence, the resulting approximation can be written as follows:
$$R({\bf A},{\bf B},x,y) = \displaystyle \frac{0.1720 + 6.6221 x^2 + 6.6221 y^2}{1 + 4.9000 x^2 + 4.9036 y^2}.$$
The maximum deviation for the corresponding problem is 0.1720 and the computational time is 725.2468 seconds ($\approx$ 12 minutes).

\begin{table}
	\centering
	\caption{Obtained values for the coefficients of $R({\bf A},{\bf B},x,y)$.}
	{\begin{tabular}{ |c|c|c| } 
			\hline
			Basis function & Coefficient & Value  \\ \hline 
			$1$ & $a_0$ & 0.1720 \\ \hline
			$x$ & $a_1$ & 0 \\ \hline
			$y$ & $a_2$ & 0  \\ \hline
			$x^2$ & $a_3$ & 6.6221  \\ \hline
			$y^2$ & $a_4$ & 6.6221  \\  \hline
			$xy$ & $a_5$ & 0 \\  \hline
			$x$ & $b_1$ & 0  \\  \hline
			$y$ & $b_2$ & 0 \\  \hline
			$x^2$ & $b_3$ & 4.9000 \\  \hline
			$y^2$ & $b_4$ & 4.9036 \\  \hline
			$xy$ & $b_5$ &  0 \\  
			\hline
	\end{tabular}}
	\label{table:ratapp}
\end{table}

In this numerical experiment, we observe the following important conclusion points. 
\begin{itemize}
	\item The multivariate polynomial approximation problem is simpler  than the multivariate generalised rational approximation and the corresponding program works much faster. It requires less than a minute to compute the polynomial approximation whereas in the case of generalised rational approximation it takes roughly 12 minutes. 
		\item Since the original function $f(x,y)$ is a nonsmooth and non-Lipschitz function, multivariate polynomial approximation is not very efficient. In particular, the shape of the multivariate generalised rational approximation is much closer to the original shape of the function: it reflects the abrupt change reasonably well compared to multivariate polynomial function.
		\item The maximum deviation of the multivariate generalised rational approximation is much smaller than the maximum deviation of the multivariate polynomial approximation. This is a promising result which indicates that the multivariate generalised rational approximation works better on nonsmooth and non-Lipschitz functions than polynomials. 
\end{itemize}

\section{Conclusions and future research directions}\label{sec:conc}
In this paper we show how a generalised rational approximation procedure, originally developed for univariate approximation,  can be extended to the case of multivariate functions. In particular, the optimisation problems remain quasiconvex and one can apply the bisection method developed for general quasiconvex problems. The main difficulty of this method is to solve the convex feasibility subproblems appearing at each step. In our computational experiments, we use the fact that in the case of approximation over a finite set of points (for example a grid), the feasibility problems can be reduced to linear programming problems.

  We also provide the results of numerical experiments that demonstrate that programs for polynomial approximation the work faster than programs for generalised rational approximation. At the same time, generalised rational approximations are better suited for approximating nonsmooth and non-Lipschitz functions.

The bisection method that we extend in this paper is simple and robust, but it is not as efficient as some other methods, in particular, the differential correction method. One of our main future research direction is to extend other numerical methods originally developed for generalised rational approximation of univariate functions to multivariate approximation. The extension of the procedure for differential correction to multivariate approximation is straightforward, but the analytical properties still need more work. 

Another important research direction is to compare methods for checking convex feasibility. This study is especially important for developing numerical methods for the case of constructing approximations in the forms of general quasilinear functions.

\section*{Acknowledgement}
This research was supported by the Australian Research Council (ARC),  Solving hard Chebyshev approximation problems through nonsmooth analysis (Discovery Project DP180100602).


\end{document}